# Power-Based Generation Expansion Planning for Flexibility Requirements

Diego A. Tejada-Arango, Germán Morales-España, Sonja Wogrin, and Efraim Centeno

*Abstract*—Flexibility requirements are becoming more relevant in power system planning due to the integration of variable Renewable Energy Sources (vRES). In order to consider these requirements Generation Expansion Planning (GEP) models have recently incorporated Unit Commitment (UC) constraints, using traditional energy-based formulations. However, recent studies have shown that energy-based UC formulations overestimate the actual flexibility of the system. Instead, power-based UC models overcome these problems by correctly modeling ramping constraints and operating reserves. This paper proposes a power-based GEP-UC model that improves the existing models. The proposed model optimizes investment decisions on vRES, Energy Storage Systems (ESS), and thermal technologies. In addition, it includes real-time flexibility requirements, and the flexibility provided by ESS, as well as other UC constraints, e.g., minimum up/down times, startup and shutdown power trajectories, network constraints. The results show that power-based model uses the installed investments more effectively than the energy-based models because it more accurately represents flexibility capabilities and system requirements. For instance, the power-based model obtains less investment (6-12%) and yet it uses more efficiently this investment because operating cost is also lower (2-8%) in a real-time validation. We also propose a semi-relaxed power-based GEP-UC model, which is at least 10 times faster than its full-integer version and without significantly losing accuracy in the results (less than 0.2% error).

*Index Terms*--generation expansion planning, unit commitment, energy storage systems, capacity expansion planning, power system planning, power generation planning.

## NOMENCLATURE

### A. Indices and sets

| | |
|---|---|
| $j \in \mathcal{J}$ | Technologies |
| $g \in \mathcal{G} \subseteq \mathcal{J}$ | Subset of thermal generation technologies |
| $v \in \mathcal{V} \subseteq \mathcal{J}$ | Subset of renewable energy sources |
| $s \in \mathcal{S} \subseteq \mathcal{J}$ | Subset of energy storage technologies |
| $b \in \mathcal{B}$ | Buses |
| $\mathcal{B}^D \subseteq \mathcal{B}$ | Subset of buses $b$ with demand consumption |
| $l \in \mathcal{L}$ | Transmission lines |
| $\omega \in \Omega$ | Scenarios |
| $k \in \mathcal{K}_g$ | Startup segments, running from 1 (the hottest) to $K_g$ (the coldest) |
| $t \in \mathcal{T}$ | time periods (e.g., hours) |

### B. Parameters

| | |
|---|---|
| $C_j^{LV}$ | Linear variable production cost [$/MWh] |
| $C_g^{NL}$ | No-load cost [$/h] |
| $C_g^{SD}$ | Shutdown cost [$] |
| $C_{gk}^{SU}$ | Startup cost for segment $k$ [$] |
| $C_g^{EM}$ | CO2 emission cost [$/MWh] |
| $C_j^{R+}, C_j^{R-}$ | Up/down reserve cost [$/MW] |
| $D_{\omega bt}^E$ | Energy demand on bus $b$ [MWh] |
| $D_{\omega bt}^P$ | Power demand on bus $b$ [MW] |
| $R_{\omega t}^+, R_{\omega t}^-$ | Up/down reserve requirement [MW] |
| $\overline{F}_l$ | Power flow limit on transmission line $l$ [MW] |
| $\overline{P}_g, \underline{P}_g$ | Maximum/minimum power output [MW] |
| $E_{gkt}^{SU}, E_{gt}^{SD}$ | Energy output during startup/shutdown [MWh] |
| $P_{gkt}^{SU}, P_{gt}^{SD}$ | Power output during startup/shutdown [MW] |
| $RU_g, RD_g$ | Ramp-up/down capability [MW/min] |
| $SU_g, SD_g$ | Startup/shutdown capability [MW] |
| $SU_g^D, SD_g^D$ | Startup/shutdown duration [h] |
| $T_{gk}^{SU}$ | Time interval limit of startup segment $k$ [h] |
| $TU_g, TD_g$ | Minimum up/down time [h] |
| $\Gamma_{lj}^J, \Gamma_{lb}$ | Shift factors for line $l$ [p.u.] |
| $EPR_s$ | Energy to power ratio [h] |
| $V_{\omega vt}^E$ | Renewable energy output profile [p.u.] |
| $V_{\omega vt}^P$ | Renewable power output profile [p.u.] |
| $\pi_\omega$ | Probability of scenario $\omega$ |
| $\overline{X}_j$ | Investment limit for technology $j$ |
| $X_j^0$ | Initial capacity for technology $j$; [# units] for $g$, and [MW] for $s$ and $v$. |

### C. Continuous non-negative variables

| | |
|---|---|
| $\hat{e}_{\omega jt}$ | Total energy output [MWh] |
| $\hat{p}_{\omega jt}$ | Total power output [MW] |
| $e_{\omega gt}$ | Energy output above minimum output [MWh] |
| $p_{\omega gt}$ | Power output above minimum output [MW] |
| $\hat{c}_{\omega st}$ | Charged energy for storage [MWh] |
| $c_{\omega st}$ | Charged power for storage [MW] |
| $r_{\omega gt}^+$ | Up capacity reserve [MW] |
| $r_{\omega gt}^-$ | Down capacity reserve [MW] |
| $\phi_{\omega st}$ | Energy storage level [MWh] |

D. A. Tejada-Arango is with the Universidad Pontificia Comillas, ICAI, Instituto de Investigación Tecnológica, Madrid, Spain (e-mail: dtejada@comillas.edu).
G. Morales-España is with the Energy Transition Studies, ECN part of TNO, Amsterdam, The Netherlands (e-mail: german.morales@tno.nl).
S. Wogrin is with the Universidad Pontificia Comillas, ICAI, Instituto de Investigación Tecnológica, Madrid, Spain (e-mail: swogrin@comillas.edu).
E. Centeno is with the Universidad Pontificia Comillas, ICAI, Instituto de Investigación Tecnológica, Madrid, Spain (e-mail: efraim@comillas.edu).

*D. Integer Variables*

| | |
|---|---|
| $u_{\omega gt}$ | Unit commitment for thermal technologies |
| $y_{\omega gt}$ | Startup for thermal technologies |
| $z_{\omega gt}$ | Shutdown for thermal technologies |
| $\delta_{\omega gkt}$ | Startup type selection for thermal technologies |
| $\gamma_{\omega st}$ | Binary decision for charging/discharging logic |
| $x_j$ | Investment decision per technology |

## II. INTRODUCTION

GENERATION Expansion Planning (GEP) is a classic long-term problem in power systems that aims at determining the optimal generation technology mix [1]. Environmental policies, such as renewable targets [2] or CO2 emission reduction [3] influence in GEP decisions, leading to the integration of vast amounts of variable Renewable Energy Sources (vRES), i.e., wind and solar, in GEP. Nevertheless, vRES integration has consequences in GEP modeling. For instance, previous studies [4]–[6] have shown the importance of including short-term dynamics on GEP decisions in order to consider the increased need of operational flexibility due to vRES integration. Therefore, correctly modeling flexibility in GEP models is crucial to reach the right conclusions in the energy transition process.

In order to consider operational flexibility in GEP, Unit Commitment (UC) modeling is needed to determine system operation [6], [7]. For example, it is known that units are being cycled more frequently due to higher vRES flexibility requirements [8]. Studies have shown that ignoring startup and shutdown processes highly overestimates the flexibility and costs of the system [9]. Another example is ramping constraints. If we focus on flexibility and want to know a good (optimal) future generation-mix and interconnection capacities for a given scenario, the GEP problem must include at least detailed ramping constraints. Moreover, operating reserve decisions have also become more relevant in GEP with the integration of vRES because they may ensure that generation technologies have an extra income to recover their investment costs through these types of ancillary services.

Despite the recent developments to consider flexibility requirements in GEP, classic GEP models are proposed using energy-based UC models. Recent studies [9]–[11] have shown that energy-based UC models cannot capture variability on demand and vRES, and even assuming that they capture it, they cannot deliver the flexibility that they promise, that is, they intrinsically and hiddenly overestimate the flexibility of the system. This is mainly because average energy levels (e.g., average level in one hour) do not provide detailed information about the instantaneous output of a generator, and constraints such as ramping-limits and demand-balance are dependent on instantaneous outputs rather than average levels. This means that more flexibility than planned by energy-based models is used in real-time operation (through operating reserves and allowing deviations on schedules) to deal with all the problems introduced by these traditional energy-based models. These problems are hidden in the formulations, and to assess really their performance, real-time simulations are required (e.g., 5-min dispatch), as it is widely discussed in [9].

More recently, power-based models have been proposed [10], [12] to overcome these problems by better exploiting the system flexibility [9], by allowing the correct modeling of ramping constraints and operating reserves [10], [11] in order to deliver the expected and actual flexibility from the generation resources. This is possible because a power-based model has a clear distinction between power and energy in its core formulation. Demand and generation are modeled as hourly piecewise-linear functions representing their instantaneous power trajectories. The schedule of a generating unit output is no longer an energy stepwise function, but a smoother piecewise power function.

Another important aspect to determine the flexibility requirements in power systems is time resolution. In order to model correctly the real operation of power systems a high resolution is needed (e.g., minutes). Current GEP models are based on hourly resolution where the underlying assumption is that it is enough to capture the variability and flexibility requirements of power systems. However, it has already been shown in [9] that real-time simulations (e.g., 5-min time step) help to determine the performance of different schedules (operational decisions) to meet the real-time flexibility requirements in the power system. This type of real-time validation is not common to be carried out because it is considered unnecessary. Nevertheless, to validate correctly flexibility capabilities and requirements of the system, this real-time evaluation is paramount [9].

In this paper, we propose a novel power-based model for GEP presenting advantages over the traditional energy-based models. The proposed model optimizes investment decisions on vRES, Energy Storage Systems (ESS), and thermal technologies. ESS are included because they represent one of the most promising options to provide flexibility in power systems in the future [13]. The main contributions of this paper are as follows:

1) We propose a power-based GEP-UC model that improves the classic energy-based models by representing more accurately the flexibility requirements of power systems (i.e., reserve decisions and ramping constraints). We present a model for ESS based on power, so it added to the power-based formulation.
2) We propose a real-time validation stage (e.g., 5-min simulation) in order to evaluate the quality of investment and operational decisions obtained with the model.
3) We also propose a semi-relaxed version of the power-based GEP-UC model, which aims to reduce the computational burden without losing accuracy in the results.
4) We show that even in the proposed semi-relaxed version, the power-based GEP-UC model obtains better performance in the real-time validation stage than the traditional energy-based models, while the investment problem is solved significantly faster.

## III. GENERATION EXPANSION MODEL FORMULATIONS

This section presents the objective function and set of constraints for the energy- and the power-based GEP-UC

formulations. These constraints include investment decisions for different generation technologies: thermal generation, ESS, and VRES. In addition, operational decisions are considered using a clustered UC formulation (i.e., aggregating similar generating units into one group or cluster), which is commonly applied in long-term planning models [7], [14], [15].

*A. Energy-Based Formulation*

The GEP seeks to minimize the investment costs plus the expected value of operating costs: production cost, up/down reserve cost, CO2 emission cost, no-load cost, shutdown cost, startup cost. Notice that $\Psi = \{x, e, \hat{e}, \hat{c}, r^+, r^-, u, y, z, \delta, \phi\}$ corresponds to the set of decision variables considered in this model.

$$\min_{\Psi} \sum_{j\in\mathcal{J}} C_j^I x_j + \sum_{\omega\in\Omega} \pi_\omega \sum_{t\in\mathcal{T}} \Big\{ \sum_{j\in\mathcal{J}} [C_j^{LV} \hat{e}_{\omega jt} + C_j^{R+} r^+_{\omega jt} + C_j^{R-} r^-_{\omega jt}] + \sum_{g\in\mathcal{G}} [C_g^{EM} \hat{e}_{\omega gt} + C_g^{NL} u_{\omega gt} + C_g^{SD} z_{\omega gt} + \sum_{k\in\mathcal{K}_g} C_{gk}^{SU} \delta_{\omega gkt}] \Big\} \quad (1)$$

The system-wide constraints are guaranteed by energy demand balance (2), transmission limits (3), and reserve requirements (4)-(5):

$$\sum_{j\in\mathcal{J}} \hat{e}_{\omega jt} - \sum_{s\in\mathcal{S}} \hat{c}_{\omega st} = \sum_{b\in\mathcal{B}^D} D^E_{\omega bt} \quad \forall \omega, t \quad (2)$$

$$-\overline{F}_l \le \sum_{j\in\mathcal{J}} \Gamma^J_{lj} \hat{e}_{\omega jt} - \sum_{s\in\mathcal{S}} \Gamma^S_{ls} c_{\omega st} - \sum_{b\in\mathcal{B}^D} \Gamma_{lb} D^E_{\omega bt} \le \overline{F}_l$$
$$\forall l, \omega, t \quad (3)$$

$$\sum_{j\in\mathcal{J}} r^+_{\omega jt} \ge R^+_{\omega t} \quad \forall \omega, t \quad (4)$$

$$\sum_{j\in\mathcal{J}} r^-_{\omega jt} \ge R^-_{\omega t} \quad \forall \omega, t \quad (5)$$

The relationship between operational and investment decisions for each technology type is guaranteed with (6) for thermal technologies, (7)-(8) for ESS, and (9) for vRES.

$$u_{\omega gt} \le X^0_g + x_g \quad \forall \omega, g, t \quad (6)$$

$$\hat{e}_{\omega st} - \hat{c}_{\omega st} + r^+_{\omega gt} \le X^0_s + x_s \quad \forall \omega, s, t \quad (7)$$

$$\hat{e}_{\omega st} - \hat{c}_{\omega st} - r^-_{\omega gt} \ge -(X^0_s + x_s) \quad \forall \omega, s, t \quad (8)$$

$$\hat{e}_{\omega vt} \le V^E_{\omega vt}(X^0_v + x_v) \quad \forall \omega, v, t \quad (9)$$

Thermal generation constraints include: commitment/ startup/ shutdown logic (10), minimum up/down times (11)-(12), startup type selection (13)-(14) (e.g., hot, warm, and cold startup), energy production limits including reserve decisions (15)-(18) (where $\mathcal{G}^1$ is defined as the thermal technologies in $\mathcal{G}$ with $TU_g = 1$), and total energy production (19). The UC formulation presented here is based on the tight and compact formulation proposed in [16]. Furthermore, Gentile et al. [17] have proven that the set of constraints (10)-(12) together with (15)-(19) is the tightest representation (i.e., convex hull) for the energy-based model.

$$u_{\omega gt} - u_{\omega g,t-1} = y_{\omega gt} - z_{\omega gt} \quad \forall \omega, g, t \quad (10)$$

$$\sum_{i=t-TU_g+1}^{t} y_{\omega gi} \le u_{\omega gt} \quad \forall \omega, g, t \in [TU_g, T] \quad (11)$$

$$\sum_{i=t-TD_g+1}^{t} z_{\omega gi} \le (X^0_g + x_g) - u_{\omega gt}$$
$$\forall \omega, g, t \in [TD_g, T] \quad (12)$$

$$\delta_{\omega gkt} \le \sum_{i=T^{SU}_{gk}}^{T^{SU}_{g,k+1}-1} z_{\omega g,t-i} \quad \forall \omega, g, k \in [1, K_g), t \quad (13)$$

$$\sum_{k\in\mathcal{K}_g} \delta_{\omega gkt} = y_{\omega gt} \quad \forall \omega, g, t \quad (14)$$

$$e_{\omega gt} + r^+_{\omega gt} \le (\overline{P}_g - \underline{P}_g) u_{\omega gt} - (\overline{P}_g - SD_g) z_{\omega g,t+1}$$
$$- \max(SD_g - SU_g, 0) y_{\omega gt} \quad \forall \omega, g \in \mathcal{G}^1, t \quad (15)$$

$$e_{\omega gt} + r^+_{\omega gt} \le (\overline{P}_g - \underline{P}_g) u_{\omega gt} - (\overline{P}_g - SU_g) y_{\omega g,t}$$
$$- \max(SU_g - SD_g, 0) z_{\omega g,t+1} \quad \forall \omega, g \in \mathcal{G}^1, t \quad (16)$$

$$e_{\omega gt} + r^+_{\omega gt} \le (\overline{P}_g - \underline{P}_g) u_{\omega gt} - (\overline{P}_g - SU_g) y_{\omega g,t}$$
$$-(\overline{P}_g - SD_g) z_{\omega g,t+1} \quad \forall \omega, g \notin \mathcal{G}^1, t \quad (17)$$

$$e_{\omega gt} - r^-_{\omega gt} \ge 0 \quad \forall \omega, g, t \quad (18)$$

$$\hat{e}_{\omega gt} = \underline{P}_g u_{\omega gt} + e_{\omega gt} \quad \forall \omega, g, t \quad (19)$$

Traditional energy-based UC formulations ignore the inherent startup (SU) and shutdown (SD) trajectories of thermal generation, assuming they start/end their production at their minimum output. Authors in [9], [10] have shown the relevance of the SU and SD processes when they are included in the scheduling optimization. Therefore, we also analyze the energy-based formulation including the SU/SD trajectories proposed in [18]. Thus, if SU/SD trajectories are considered then (19) is replaced by (20).

$$\hat{e}_{\omega gt} = \underbrace{\sum_{k=1}^{K_g} \sum_{i=1}^{SU^D_{gk}} E^{SU}_{gki} \delta_{\omega gk,(t-i+SU^D_{gk}+1)}}_{\text{Startup trajectory}} +$$
$$\underbrace{\sum_{i=1}^{SD^D_g} E^{SD}_{gi} z_{\omega g,(t-i+1)}}_{\text{Shutdown trajectory}} + \underbrace{\underline{P}_g u_{\omega gt} + e_{\omega gt}}_{\text{Output when being up}} \quad \forall \omega, g, t \quad (20)$$

ESS constraints include: logic to avoid charging and discharging at the same time (21)-(22), the definition of the storage inventory level (23), storage limits including reserve (24)-(25).

$$\hat{c}_{\omega st} \le (1 - \gamma_{\omega st}) \cdot (X^0_s + \overline{X}_s) \quad \forall \omega, s, t \quad (21)$$

$$\hat{e}_{\omega st} \le \gamma_{\omega st} \cdot (X^0_s + \overline{X}_s) \quad \forall \omega, s, t \quad (22)$$

$$\phi_{\omega st} = \phi_{\omega s,t-1} + \eta_s \hat{c}_{\omega st} - \hat{e}_{\omega st} \quad \forall \omega, s, t \quad (23)$$

$$\phi_{\omega st} \le EPR_s(X^0_s + x_s) - \sum_{i=t-1}^{t} r^-_{\omega gi} \quad \forall \omega, s, t \quad (24)$$

$$\phi_{\omega st} \ge \sum_{i=t-1}^{t} r^+_{\omega gi} \quad \forall \omega, s, t \quad (25)$$

Flexibility requirements in the power system are represented by ramping constraints including reserve decisions. In order to guarantee that scheduled reserves are feasible to provide at $\tau$-min (e.g., $\tau$=5) using the energy-based formulation, it is necessary to consider the ramping capability at $\tau$-min. For instance, ramp capability limits imposed with (26)-(27) consider the reserve that thermal technologies can provide at $\tau$-min. ESS ramp capability limits (28)-(29) consider the charged energy in addition to the energy output (i.e., discharged energy). Notice that (28)-(29) allow ESS to switch from charging to discharging within the ramp limit, as in [19].

$$(e_{\omega gt} + r^+_{\omega gt}) - e_{\omega g,t-1} \le \tau RU_g u_{\omega gt} \quad \forall \omega, g, t \quad (26)$$

$$(e_{\omega gt} - r^-_{\omega gt}) - e_{\omega g,t-1} \ge -\tau RD_g u_{\omega g,t-1} \quad \forall \omega, g, t \quad (27)$$

$$(\hat{c}_{\omega st} - \hat{c}_{\omega s,t-1}) + (\hat{e}_{\omega st} - \hat{e}_{\omega s,t-1}) + r^+_{\omega st} \le$$
$$\tau RU_g(X^0_s + x_s) \quad \forall \omega, s, t \quad (28)$$

$$(\hat{c}_{\omega st} - \hat{c}_{\omega s,t-1}) + (\hat{e}_{\omega st} - \hat{e}_{\omega s,t-1}) - r^-_{\omega st} \ge$$
$$-\tau RD_g(X^0_s + x_s) \quad \forall \omega, s, t \quad (29)$$

*B. Power-Based Formulation*

This section shows the GEP-UC equations in terms of power. However, some of the terms in these equations are naturally linked to energy. For instance, the objective function (30) considers the so-called calculated energy $\hat{e}_{\omega jt}$ to determine the variable cost and CO2 emission cost. The calculated energy is determined using the power output variables $\hat{p}_{\omega jt}$ in (31). In addition, for ESS the charged energy $\hat{c}_{\omega st}$ is also determined using the charged power in (32). Notice that $\Lambda = \{x, p, \hat{p}, e, c, \hat{c}, r^+, r^-, u, y, z, \delta, \phi\}$ corresponds to the set of decision variables in this model.

$$\min_{\Lambda} \sum_{j \in \mathcal{J}} C_j^I x_j + \sum_{\omega \in \Omega} \pi_\omega \sum_{t \in \mathcal{T}} \Big\{ \sum_{j \in \mathcal{J}} [C_j^{LV} \hat{e}_{\omega jt} + C_j^{R+} r_{\omega jt}^+ + C_j^{R-} r_{\omega jt}^-] + \sum_{g \in \mathcal{G}} [C_g^{EM} \hat{e}_{\omega gt} + C_g^{NL} u_{\omega gt} + C_g^{SD} z_{\omega gt} + \sum_{k \in \mathcal{K}_g} C_{gk}^{SU} \delta_{\omega gkt}] \Big\} \quad (30)$$

$$\hat{e}_{\omega jt} = \frac{\hat{p}_{\omega jt} + \hat{p}_{\omega j,t-1}}{2} \quad \forall \omega, j, t \quad (31)$$

$$\hat{c}_{\omega st} = \frac{c_{\omega st} + c_{\omega s,t-1}}{2} \quad \forall \omega, s, t \quad (32)$$

Demand balance constraint (33) and power-flow transmission limits (34) also use the power output instead of energy output. Reserve requirements (4)-(5) remain the same because they are already expressed in terms of power.

$$\sum_{j \in \mathcal{J}} \hat{p}_{\omega jt} - \sum_{s \in \mathcal{S}} c_{\omega st} = \sum_{b \in \mathcal{B}^D} D_{\omega bt}^E \quad \forall \omega, t \quad (33)$$

$$-\overline{F}_l \leq \sum_{j \in \mathcal{J}} \Gamma_{lj}^J \hat{p}_{\omega jt} - \sum_{s \in \mathcal{S}} \Gamma_{ls}^S c_{\omega st} - \sum_{b \in \mathcal{B}^D} \Gamma_{lb} D_{\omega bt}^E \leq \overline{F}_l \quad \forall l, \omega, t \quad (34)$$

In terms of the relationship between operational and investment decisions, thermal unit constraint (6) remains the same. However, constraints for ESS and vRES technologies change to (35)-(36) and (37), respectively.

$$\hat{p}_{\omega st} - c_{\omega st} + r_{\omega gt}^+ \leq X_s^0 + x_s \quad \forall \omega, s, t \quad (35)$$

$$\hat{p}_{\omega st} - c_{\omega st} - r_{\omega gt}^- \geq -(X_s^0 + x_s) \quad \forall \omega, s, t \quad (36)$$

$$\hat{p}_{\omega vt} \leq V_{\omega vt}^P (X_v^0 + x_v) \quad \forall \omega, v, t \quad (37)$$

Unit commitment constraints (10)-(14) do not change in the power-based formulation. Equations (38)-(39) limit the power output of thermal technologies. The total power output constraint is different depending whether it is a quick- or slow-start unit. Quick-start technologies $\mathcal{G}^F$ are thermal generators that can startup/shutdown within one hour (i.e., $SU_{gk}^D = SD_g^D \leq 1$), while slow-start technologies $\mathcal{G}^S$ are those with a SU/SD duration greater than one hour as well as a SU/SD capacity equal to the minimum power output (i.e., $SU_g = SD_g = \underline{P}_g$). Therefore, the total power output of slow-start technologies considers SU/SD trajectories (41), whereas (40) for quick-start technologies does not. For a better understanding of the modeling of quick- and slow-start technologies, the reader is referred to [12], [17]. The formulation presented here is based on the tight and compact formulation proposed in [10]. Furthermore, Morales-España et al. [12] has proven that the set of constraints (10)-(12) together with (38)-(41) is the tightest possible representation (i.e., convex hull) for the power-based model.

$$p_{\omega gt} + r_{\omega gt}^+ \leq (\overline{P}_g - \underline{P}_g) u_{\omega gt} - (\overline{P}_g - SD_g) z_{\omega g,t+1} + (SU_g - \underline{P}_g) y_{\omega g,t+1} \quad \forall \omega, g, t \quad (38)$$

$$p_{\omega gt} - r_{\omega gt}^- \geq 0 \quad \forall \omega, g, t \quad (39)$$

$$\hat{p}_{\omega gt} = \underline{P}_g (u_{\omega gt} + y_{\omega g,t+1}) + p_{\omega gt} \quad \forall \omega, g \in \mathcal{G}^F, t \quad (40)$$

$$\hat{p}_{\omega gt} = \underbrace{\sum_{k=1}^{K_g} \sum_{i=1}^{SU_{gk}^D} P_{gki}^{SU} \delta_{\omega gk,(t-i+SU_{gk}^D+2)}}_{\text{Startup trajectory}} + \underbrace{\sum_{i=2}^{SD_g^D+1} P_{gi}^{SD} z_{\omega g,(t-i+2)}}_{\text{Shutdown trajectory}} + \underbrace{\underline{P}_g (u_{\omega gt} + y_{\omega g,t+1}) + p_{\omega gt}}_{\text{Output when being up}} \quad \forall \omega, g \in \mathcal{G}^S, t \quad (41)$$

ESS constraints for storage level (23) and storage level limits including reserve (24)-(25) continue the same. Nevertheless, the logic to avoid charging and discharging at the same time (42)-(43) is updated to consider the power output and charged power.

$$c_{\omega st} \leq (1 - \gamma_{\omega st}) \cdot (X_s^0 + \overline{X}_s) \quad \forall \omega, s, t \quad (42)$$

$$\hat{p}_{\omega st} \leq \gamma_{\omega st} \cdot (X_s^0 + \overline{X}_s) \quad \forall \omega, s, t \quad (43)$$

One of the main advantages of power-based formulation is that it allows to describe a more detailed set of constraints to represent the flexibility requirements, which are described in terms of power instead of energy. The proposed power-based equations in [10] ensure that reserves can be provided at any time within the hour by guaranteeing that the reserve does not exceed the ramp-capability at $\tau$-min (e.g., $\tau$=5 min) and power-capacity limits at the end of the hour (i.e., 60 min). Therefore, (44)-(45) guarantee that $\tau$-min ramp capability is ensured for thermal technologies, while (46)-(47) guarantee the power-capacity limit for both $\tau$-min and at the end of the hour. For a detailed explanation of how the reserve constrains are obtained, the reader is referred to [10]. Although [10] shows the case for thermal technologies, we use the same concepts and extend the concept for ESS in (48)-(51).

$$\frac{\tau(p_{\omega gt} - p_{\omega g,t-1})}{60} + r_{\omega gt}^+ \leq \tau RU_g u_{\omega gt} \quad \forall \omega, g, t \quad (44)$$

$$\frac{\tau(p_{\omega gt} - p_{\omega g,t-1})}{60} - r_{\omega gt}^- \geq -\tau RD_g u_{\omega g,t-1} \quad \forall \omega, g, t \quad (45)$$

$$\frac{\tau p_{\omega gt} + (60-\tau) p_{\omega g,t-1}}{60} + r_{\omega gt}^+ \leq (\overline{P}_g - \underline{P}_g) u_{\omega gt} \quad \forall \omega, g, t \quad (46)$$

$$\frac{\tau p_{\omega gt} + (60-\tau) p_{\omega g,t-1}}{60} - r_{\omega gt}^- \geq 0 \quad \forall \omega, g, t \quad (47)$$

$$\frac{\tau(c_{\omega st} - c_{\omega s,t-1})}{60} + \frac{\tau(\hat{p}_{\omega st} - \hat{p}_{\omega s,t-1})}{60} + r_{\omega st}^+ \leq \tau RU_g (X_s^0 + x_s) \quad \forall \omega, s, t \quad (48)$$

$$\frac{\tau(c_{\omega st} + \hat{p}_{\omega st}) + (60-\tau)(c_{\omega s,t-1} + \hat{p}_{\omega s,t-1})}{60} + r_{\omega st}^+ \leq X_s^0 + x_s \quad \forall \omega, s, t \quad (49)$$

$$\frac{\tau(c_{\omega st} - c_{\omega s,t-1})}{60} + \frac{\tau(\hat{p}_{\omega st} - \hat{p}_{\omega s,t-1})}{60} - r_{\omega st}^- \geq -\tau RD_g (X_s^0 + x_s) \quad \forall \omega, s, t \quad (50)$$

$$\frac{\tau(c_{\omega st} + \hat{p}_{\omega st}) + (60-\tau)(c_{\omega s,t-1} + \hat{p}_{\omega s,t-1})}{60} - r_{\omega st}^- \geq 0 \quad \forall \omega, s, t \quad (51)$$

IV. SYSTEM FLEXIBILITY EVALUATION

As mentioned in the previous section, two main formulations are analyzed for GEP: the traditional energy-based (EB), and the power-based formulation (PB). We also analyze the traditional energy-based using SU/SD trajectories (EBs). Table I shows a summary with all the equations that define these models. All models include an hourly UC (either energy- or power-based) in order to consider operating constraints, involving those related to the power system flexibility (i.e., ramping and reserve constraints). In order to measure the quality of the obtained solution under real-time flexibility requirements, we carry out an evaluation of investment and operational decisions through a simulation using the same scenarios as in the GEP-UC hourly optimization (in-sample simulation). This evaluation allows us to establish the problems associated to each formulation rather than those associated to the uncertainty representation by itself. The complete procedure to calculate investment decisions and ex-post real-time evaluation is shown in Fig. 1 (top). During stage 1, the investment and hourly UC schedule are optimized solving the formulations shown in Section III. Then, investment, commitment, and reserve decisions are fixed. Stage 2 tests the results through a real-time simulation model, using a 5-min optimal dispatch (emulating real-time markets as in [9]) in order to evaluate the GEP-UC solution. Dispatch decisions (e.g., production, charge/discharge) obtained in stage 2 are called redispatches, allowing us to evaluate the deviations with respect to the stage 1. This is called the *integer approach*. In addition, we proposed a *semi-relaxed approach* for the power-based formulation, which is shown in Fig. 1 (bottom). Here we split stage 1 in two. First, the stage 1a solves the power-based formulation considering integer investment decisions and continuous UC decisions ($u_{\omega gt}, y_{\omega gt}, z_{\omega gt}$). This approximation allows to solve the GEP problem much faster.

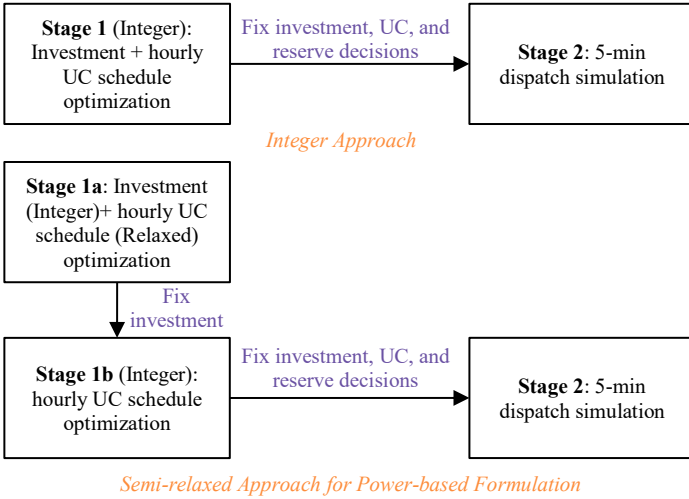

Fig. 1. Stage sequence for integer (top) and semi-relaxed (bottom) approaches.

TABLE I
GEP-UC Models

| Equations | EB | EBs | PB | SR-PB |
|---|---|---|---|---|
| Objective function | (1) | | (30) | |
| System constraints | (2)-(5) | | (4)-(5), (31)-(34) | |
| Investment constraints | (6),(7)-(8),(9) | | (6), (35)-(36),(37) | |
| UC constraints | (10)-(14) | | | |
| Thermal unit constraints | (15)-(18) | | (38)-(39) | |
| Total output thermal technologies | (19) | (20) | (40)-(41) | |
| ESS constraints | (21)-(25) | | (23)-(25),(43)-(42) | |
| Constraints for flexibility requirements ($\tau = 5$min) | (26)-(29) | | (44)-(51) | |
| Integer variables | $u_{\omega gt}, y_{\omega gt}, z_{\omega gt}, \gamma_{\omega st}, \delta_{\omega gkt}, x_j$ | | | Stage 1a: $x_j$ Stage 1b: $u_{\omega gt}, y_{\omega gt}, z_{\omega gt}, \gamma_{\omega st}, \delta_{\omega gkt}$ |

Then investment decisions are fixed in stage 1b, where the power-based formulation is solved considering integer UC decisions. Once again, investment, unit commitment, and reserve decisions are fixed to simulate a 5-min optimal dispatch.

## V. CASE STUDIES

To evaluate the performance of the different approaches, we use two case studies: a modified IEEE 118-bus test system and a stylized Dutch power system in target year 2040. Input data for both case studies is available online at [20], including the 5-min demand and renewable production profiles. Both case studies are solved considering a green-field investment approach (i.e., no initial capacity) for thermal generation and ESS investment, while the vRES capacity is predefined.

The modified IEEE 118-bus test system is described in Morales-España [21] for a time span of 24 h. This system was originally conceived for UC problems and it has 118 buses, 186 transmission lines, 91 loads, 54 slow-start thermal technologies, 10 quick-start technologies, and three buses with wind production. Nevertheless, we adapt this case study for GEP problems. Thermal unit investments are allowed in buses where there was a unit connected in the initial UC problem. In addition, ESS investment decisions are available in three types of technologies (PSH, CAES, and Li-ION) for buses with renewable production. The total (5-min) load average is 3578.6MW, it has a peak of 5117.5MW and a minimum of 1435.4MW.

The stylized Dutch system case study for year 2040 is mainly based on the information available in the *Ten Year Network Development Plan 2018* [22] (e.g., hourly demand profile, renewable capacity, technical characteristics and available technologies). However, the wind and solar profiles were taken from [23], [24] since this information is not available in [22]. Instead of solving 8760 h for the whole year, we have selected four representative weeks using the proposed method in [25] and k-medoids clustering technique [26]. Other authors [27], [28] have proposed different approaches to select the representative periods (e.g., weeks or days) that are compatible with the proposed GEP-UC models in this paper. Each representative week is considered as one scenario in the optimization problem, and the scenario probability is obtained from the clustering process. For investment decisions, four different thermal generation technologies are considered, Combined Heat and Power (CHP), combined cycle gas turbine (CCGT), open cycle gas turbine (OCGT), and Light Oil (Oil). Moreover, three ESS (PSH, CAES, Li-ION) technologies are considered for investment decisions.

For each case study, four different models are implemented: traditional energy-based (EB), energy-based including SU/SD power trajectories (EBs), the proposed power-based formulation (PB), and the semi-relaxed power-based formulation (SR-PB). Table I shows the summary with all the implemented models. All models consider $\tau = 5$min for constraints associated to flexibility constraints.

All optimizations were carried out using Gurobi 8.1 on an Intel®-Core™ i7-4770 (64-bit) 3.4-GHz personal computer with 16GB of RAM memory. The problems are solved until they reach an optimality tolerance of 0.1%.

## VI. Results

### A. Modified IEEE 118-bus System

Table II shows the main results for each model. The total investment cost (ESS + Thermal) is higher in the classic EB model than the one obtained with the PB model. Generally, increasing the investments lowers operating cost. Nevertheless, here we obtain a counterintuitive result. Even though the classic EB model invests more (6%), the operating cost is worse than the one in the PB model (15%). Moreover, the $CO_2$ emissions and curtailment are also higher in the classic EB model, despite its higher capacity in clean ESS and lower capacity in thermal technologies. This is also a counterintuitive result, because at a first glance, less thermal generation should pollute less, and more storage should allocate more renewables. However, this result is related to how the technology mix is selected in each model. Therefore, it is not only a matter of how much the model invests, it is also a matter of how the technology mix is selected, see Table III. For instance, although the total coal capacity is higher in the proposed PB model, the actual total coal production is lower (7%) than the one in the classic EB model, see Table IV. This is compensated by a higher use of wind, gas (that have less $CO_2$ emission factor) and oil, which overall results in less $CO_2$ emissions. As mentioned in Section III. B. the PB model equations allow to schedule the thermal technologies in a way that correctly represents the requirements and actual availability of system's flexibility, such as the load ramps. The results show the benefits of accurately considering the flexibility requirements and of correctly modelling the flexibility capabilities of the system by modelling in terms of power instead of energy.

The EBs model improves the classic EB model by including the SU/SD power-based ramps. In stage 1, the total cost in the EBs model is 8.5% lower than the classic EB model. However, it is still 4% higher than the PB model and with more curtailment (5.7 times). The EBs technology mix is also different, as it invests more in PHS and coal (Table III). And yet, the PB model allocates more wind with less ESS, see Table IV. Therefore, the PB model invests more efficiently due to the more accurate representation of flexibility requirements and capabilities of the power system.

Regarding the CPU time, the PB model is faster than its energy counterparts (2.4 and 1.5 times respectively). Nevertheless, for large-scale investment decision problems, the integer nature of the UC variables especially could make the problem intractable to solve. Therefore, the proposed SR-PB models aims at overcoming this situation. For instance, it solves the problem 9 times faster than the PB model and with only a 0.2% difference in the objective function. Moreover, the difference in the $CO_2$ emissions is only 0.4%. The main difference appears in the curtailment (90%) due to the increase in the investment made by the SR-PB that allows to reduce the operating cost by increasing wind production. When the SR-PB and the EB are compared, it may be concluded that the even the semi-relaxed version of the power-based model (i.e., SR-PB) shows better performance than the discrete version of the energy-based models (i.e., EB and EBs). In other words, the SR-PB model has a lower total cost than the EB model, investing and operating with lower cost, while simultaneously solving 21+ times faster.

TABLE II
IEEE 118-bus System: Performance for each formulation

| | Result | EB | EBs | PB | SR-PB |
|---|---|---|---|---|---|
| Stage 1 | Total Cost [M$] | 10.15 | 9.29 | 8.94 | 8.96[†] |
| | ESS Invest Cost [M$] | 0.43 | 0.35 | 0.19 | 0.17 |
| | Therm. Invest Cost [M$] | 1.01 | 1.42 | 1.17 | 1.24 |
| | Operating Cost [M$] | 8.71 | 7.52 | 7.58 | 7.55[†] |
| | CO2 emissions [ton] | 63.11 | 53.06 | 53.98 | 53.74 |
| | Curtailment [%] | 5.76 | 4.18 | 0.73 | 0.70 |
| | CPU Time [s] | 10717 | 6767 | 4478 | 500 |
| Stage 2 | Operating Cost [M$] | 8.22 | 7.53 | 7.58 | 7.55 |
| | Total Cost [M$] | 9.66 | 9.30 | 8.94 | 8.96 |
| | CO2 emissions [ton] | 59.31 | 52.48 | 53.95 | 53.71 |
| | Curtailment [%] | 0.00 | 0.00 | 0.60 | 0.62 |

[†] Values from Stage 1b

TABLE III
Technology investment decisions [MW]

| Technology | EB | EBs | PB | SR-PB |
|---|---|---|---|---|
| PSH | 1250 | 1000 | 500 | 441 |
| CAES | 0 | 0 | 0 | 0 |
| Li-ION | 150 | 150 | 150 | 150 |
| GAS | 360 | 600 | 420 | 480 |
| COAL | 4380 | 6080 | 5030 | 5330 |
| OIL | 50 | 100 | 100 | 100 |

TABLE IV
Technology production decisions [MWh]

| Technology | EB | EBs | PB | SR-PB |
|---|---|---|---|---|
| PSH | 7352 | 5944 | 2449 | 2019 |
| CAES | 0 | 0 | 0 | 0 |
| Li-ION | 1053 | 1003 | 1035 | 1033 |
| GAS | 494 | 2719 | 2482 | 2680 |
| COAL | 67540 | 63939 | 62913 | 62570 |
| OIL | 52 | 900 | 950 | 900 |
| WIND | 18880 | 19196 | 19887 | 20018 |

The results in Table II for the stage 2 are also showing interesting information: comparing the operating cost between stage 1 and 2, the classic EB shows a decrease of 6%, while in the other models remain almost the same. Moreover, the curtailment is also reduced from stage 1 to stage 2 in both energy-based models, while it remains almost the same in the power-based models. These results suggest that the obtained schedule in stage 1 with energy-based models leads to more redispatches in the technologies in stage 2. Fig. 2 illustrates this situation with the deviation with respect to the hourly thermal production obtained in stage 2 for each model. In both energy-based models, downward deviations are higher than upward deviations, which explains why the operating cost is reduced from stage 1 to stage 2 in the classic EB model as well as the reduction on the curtailment for both energy-based models. The power-based models show deviations in both directions lower than 3%, which means that the hourly schedule (stage 1) is better fitted for the 5-min real-time operation (stage 2). This high deviation of the energy-based models is due to its intrinsic incapability to accurately represent the flexibility needs and capabilities. These conclusions are aligned with those in [9] where different case studies where carried out disregarding investment decisions.

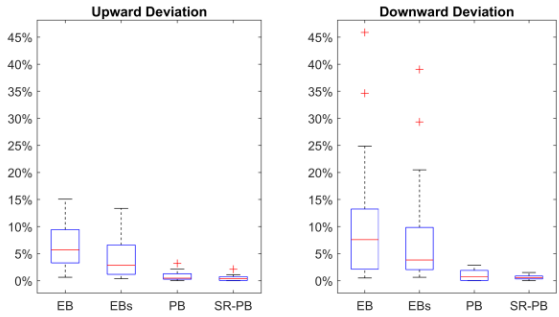

Fig. 2. Stage 2 deviation in scheduled thermal output.

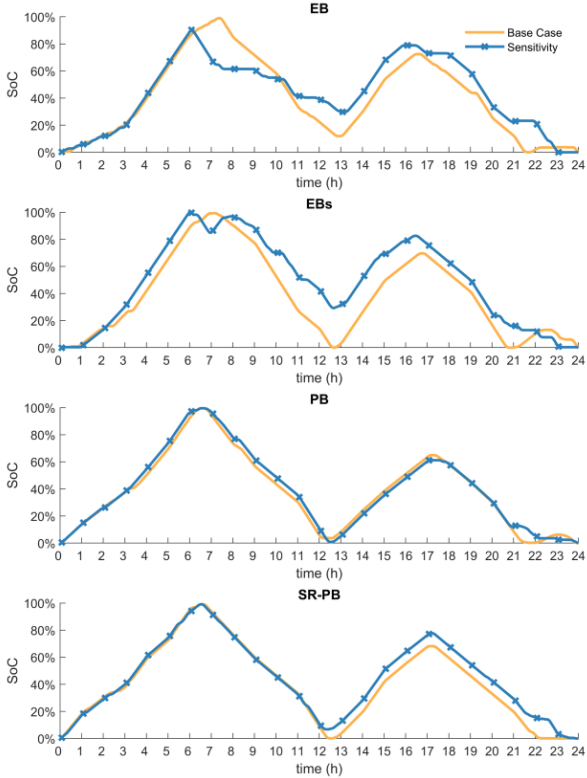

Fig. 3. Battery SoC in Stage 2 obtained for each model.

TABLE V
IEEE 118-bus System: Stage 2 – sensitivity results

| | Result | EB | EBs | PB | SR-PB |
|---|---|---|---|---|---|
| Stage 2 | Operating Cost [M$] | 8.35 | 7.66 | 7.60 | 7.55 |
| | Total Cost [M$] | 9.79 | 9.43 | 8.96 | 8.96 |
| | CO2 emissions [ton] | 59.89 | 52.71 | 54.04 | 53.73 |
| | Curtailment [%] | 1.99 | 0.98 | 0.62 | 0.13 |

Notice that ESS plays an important role in the reschedules made in stage 2. Therefore, we run a sensitivity case in which the State-of-Charge (SoC) at the end of each hour is a lower bound for the ESS in the stage 2. This limits the reschedules made in this stage, increasing the operating cost. Table V shows that situation, where with this additional constraint the operating cost, CO2 emissions and curtailment are higher than in the base case. It is important to highlight that in this sensitivity case energy-type models cannot reduce the curtailment to zero as it was in the base case. Therefore, the flexibility provided by the ESS was partly responsible for the reduction of the curtailment between stage 1 and 2 in this type of models. Fig. 3 shows the SoC in the batteries during stage 2 for the base case and the sensitivity case.

TABLE VI
Stylized Dutch System: Performance for each formulation

| | Result | EB | EBs | PB | SR-PB |
|---|---|---|---|---|---|
| Stage 1 | Total Cost [M$] | 73.18 | 70.39 | 68.14 | 68.16† |
| | ESS Invest Cost [M$] | 13.47 | 11.15 | 10.53 | 10.88 |
| | Therm. Invest Cost [M$] | 13.79 | 14.12 | 13.43 | 13.47 |
| | Operating Cost [M$] | 45.92 | 45.12 | 44.18 | 43.81† |
| | CO2 emissions [ton] | 112.10 | 98.06 | 89.46 | 88.77 |
| | Curtailment [%] | 44.72 | 45.47 | 45.46 | 45.34 |
| | CPU Time [s] | 571 | 161 | 131 | 60 |
| Stage 2 | Operating Cost [M$] | 45.76 | 46.61 | 44.90 | 44.52 |
| | Total Cost [M$] | 73.02 | 71.88 | 68.86 | 68.87 |
| | CO2 emissions [ton] | 107.73 | 100.01 | 94.29 | 93.44 |
| | Curtailment [%] | 47.88 | 48.35 | 48.34 | 45.39 |

† Values from Stage 1b

TABLE VII
Stylized Dutch System: Sensitivity to Ramp Capacity

| | Result | EB | EBs | PB | SR-PB |
|---|---|---|---|---|---|
| Stage 1 | Total Cost [M$] | 70.51 | 67.93 | 67.60 | 67.61† |
| | ESS Invest Cost [M$] | 13.35 | 10.97 | 10.66 | 10.74 |
| | Therm. Invest Cost [M$] | 13.47 | 13.47 | 13.43 | 13.47 |
| | Operating Cost [M$] | 43.69 | 43.49 | 43.51 | 43.40† |
| | CO2 emissions [ton] | 103.10 | 90.84 | 88.46 | 88.17 |
| | Curtailment [%] | 44.32 | 45.22 | 45.16 | 45.19 |
| | CPU Time [s] | 142 | 130 | 100 | 43 |
| Stage 2 | Operating Cost [M$] | 44.37 | 45.92 | 44.25 | 44.21 |
| | Total Cost [M$] | 71.19 | 70.36 | 68.34 | 68.42 |
| | CO2 emissions [ton] | 100.76 | 94.35 | 93.07 | 93.07 |
| | Curtailment [%] | 44.62 | 45.30 | 45.24 | 45.28 |

† Values from Stage 1b

The difference between both results in each model shows how the energy-type models were taking advantage of the ESS to reduce the operating cost in stage 2 at the cost of more rescheduling in the thermal technologies.

### B. Stylized Dutch System

Table VI shows the results for a stylized Dutch power system. The main conclusions drawn from the previous case study remain valid. That is, the classic EB model obtains the most expensive investment, and the operating cost is also the highest, while also resulting in the highest CO2 emissions. The amount of ESS invested in the EB model is also the highest, hence allowing it to obtain less curtailment than PB in the stage 2. Nevertheless, still the PB model results in the lowest total cost in both stages and solves the GEP problem faster than EB. In addition, the SR-PB further reduces the CPU time without losing accuracy in the results. Therefore, modeling flexibility requirements with the PB model leads to a better solution than the classic EB model.

In addition to the base case shown in Table VI, Table VII shows a sensitivity where ramp capabilities of thermal technologies are twice than before, i.e., thermal technologies are now much more flexible. As the flexibility of the thermal resources increases, the difference between energy-based and power-based models decreases. For instance, the difference between the EB the PB models changes from 7.4% to 4.3%. Therefore, if the power system does not have ramp problems, i.e., flexibility is not a problem in general, the difference between energy-based and power-based models is less significant. However, if flexibility is a limited resource and needs to be correctly managed, then the power-based models are the right option to obtain the capacity expansion planning for the system.

## VII. CONCLUSIONS

This paper proposes a power-based model to determine the GEP, including energy storage technologies. The proposed power-based model uses the installed investments more efficiently and more effectively as 1) it represents the reality of flexibility requirements of the power system more adequately, and 2) it adequately exploits the flexibility capabilities of the system. That is, the decisions made with the power-based model simultaneously yield lower investment costs, operating cost, $CO_2$ emissions, and renewable curtailment with respect to the energy-based model. This is mainly because the energy-based model overestimates flexibility capabilities, failing to capture the flexibility requirements such as load and vRES ramps even in a deterministic approach (i.e., without uncertainty on demand, or renewable production). Moreover, the advantages of the power-based approach could become much more significant considering uncertainty [9]. Therefore, correctly modeling the system flexibility changes the optimal expansion capacity decisions. For instance, the power-based model obtains less total investment (6-12%) because it is more accurate in the representation of ramping characteristics for generation resources (e.g., thermal technologies and ESS), which leads to less operating cost (2-8%) in the real-time validation. In addition, the power-based model has computational advantages in terms of CPU time. The results show that the power-based model is 2 to 4 times faster than the energy-based model. We also have demonstrated that the semi-relaxed power-based model is even faster (10 to 21 times) without losing accuracy in the results compared with the non-relaxed power-based model (less than 0.2% objective function error). This is relevant for applications with large-scale long-term capacity expansion planning problems where relaxed models are more often used due to computational power limitations.

The results show an important insight for ISOs because, even without uncertainty, the current energy-based models impose more rescheduling in the real-time operation than the power-based models. For planning authorities this is also important because decisions made with power-based models lead to a generation technology mix that is better adapted to real-time system operation.